\documentclass[11pt]{article}

\usepackage{amsmath}
\usepackage{amsthm}
\usepackage{amssymb}
\usepackage{amscd}
\usepackage[all]{xy}

\title{An example of Fourier--Mukai partners of minimal elliptic surfaces 
\footnote{2000 Mathematics Subject Classification. 
Primary: 14J27, 14E26}} 
\author{Hokuto Uehara}
\date{\empty}

\theoremstyle{plain}
\newtheorem{prop}{Proposition}[section]

\newtheorem{thm}[prop]{Theorem}

\newtheorem{cla}[prop]{Claim}

\newtheorem{mthm}{Main Theorem}

\theoremstyle{definition}

\theoremstyle{remark}

\DeclareFontFamily{U}{UWCyr}{}
\DeclareFontShape{U}{UWCyr}{m}{n}{
   <5> <6> <7> <8> <9> gen * wncyr
   <10> <10.95> <12> <14.4> <17.28> <20.74> <24.88> wncyr10
   }{}
\DeclareMathAlphabet\cy{U}{UWCyr}{m}{n}

\newcommand{\Aut}{\operatorname{Aut}}
\newcommand{\Spec}{\operatorname{Spec}}
\newcommand{\Gal}{\operatorname{Gal}}
\newcommand{\PP}{\mathbb P}
\newcommand{\Z}{\mathbb Z}
\newcommand{\Q}{\mathbb Q}
\newcommand{\C}{\mathbb C}

\begin{document}

\maketitle
\begin{abstract}\noindent 
Let $X$ and $Y$ be smooth projective varieties over $\C$. We say that $X$ and $Y$ are \emph{D-equivalent} 
(or, $X$ is a \emph{Fourier--Mukai partner} of $Y$) if their derived categories of bounded complexes of 
coherent sheaves are equivalent as triangulated categories.
The aim of this short note is to find an example of mutually D-equivalent but not isomorphic relatively 
minimal elliptic surfaces.
\end{abstract}

\section{Introduction}
Let $X$ be a smooth projective variety over $\C$. The derived category $D(X)$ of $X$ is a triangulated category
whose objects are bounded complexes of coherent sheaves on $X$.
A \emph{Fourier--Mukai (FM) transform} relating smooth projective varieties $X$ and $Y$ is an equivalence
of triangulated categories $\Phi: D(X)\to D(Y)$.
If there exists an FM transform relating $X$ and $Y$, we call $X$ an \emph{FM partner} of $Y$.
We also say that $X$ and $Y$ are \emph{D-equivalent}. Moreover we say that $X$ and $Y$ are \emph{K-equivalent}
if there exist a smooth projective variety $Z$ and birational morphisms $f:Z\to X$, $g:Z\to Y$
such that $f^* K_X\sim g^*K_Y$. It is conjectured by Kawamata (Conjecture 1.2 in \cite{Ka02}) 
that given birationally equivalent smooth projective varieties $X$ and $Y$, 
they are D-equivalent if and only if they are K-equivalent. 
In this note, we construct a counterexample to his conjecture. More precisely,
we have:
\begin{mthm}\label{main result}
\begin{enumerate}\renewcommand{\labelenumi}{(\roman{enumi})}
\item
Let $p$ be a positive integer. Then there is a rational elliptic surface
$S(p)$ such that $S(p)$ has a singular fiber of type $_{p}I_{0}$ and at least
three non-multiple singular fibers of different Kodaira's types.
\item
Let $N$ be a positive integer and $p$ a prime number such that
$p > 6(N-1)+1$. Then there are rational elliptic surfaces $T_i$, 
$(1\le i\le N)$ such that $T_i \not\cong T_j$ for $i\ne j$ and 
every $T_i$ is an FM partner of $S(p)$. As a special case,
$S = S(11)$ has an FM partner $T$ such that $T \not\cong S$. These $S$ and
$T$ are birational, D-equivalent but not K-equivalent.
\end{enumerate}
\end{mthm}
\noindent
Note that if $X$ and $Y$ are K-equivalent, they are isomorphic in codimension $1$ (Lemma 4.2, \cite{Ka02}).
In particular, if surfaces $S$ and $T$ are not isomorphic, they are not K-equivalent. Hence, in (ii), 
the statement for $p=11$ follows from the one for arbitrary $p$. 

Before ending Introduction, we give a few remarks to Main Theorem.
For a smooth projective variety $X$, it is an interesting problem to find the set of 
isomorphic classes of FM partners of $X$. In connection with this problem, we have the following.
\begin{thm}[Theorem 1.1, \cite{BM01} and Theorem 1.6, \cite{Ka02}]\label{BMK}
Assume that $X$ and $Y$ are D-equivalent smooth projective surfaces
but not isomorphic to each other.
Then we know that one of the following holds.
\begin{enumerate}\renewcommand{\labelenumi}{(\roman{enumi})}
\item $X$ and $Y$ are K3 surfaces.
\item $X$ and $Y$ are abelian surfaces.
\item $X$ and $Y$ are elliptic surfaces with the non-zero Kodaira dimension $\kappa (X)=\kappa (Y)$.
\end{enumerate}
\end{thm}
\noindent
Using Theorem \ref{BMK}, we obtain the complete answer to the problem mentioned above in dimension $2$
(\cite{BM01}, see also \cite{Ka02}).
It is well-known that the cases (i) and (ii) in Theorem \ref{BMK} really occur. More strongly, we have: 
\begin{thm}[\cite{Oguiso02} and \cite{HLOY}]\label{Oguiso}
Let $N$ be a positive integer. Then there are K3 (respectively, abelian) surfaces $T_i$, $(1\le i\le N)$ 
such that $T_i \not\cong T_j$ for $i\ne j$ and all $T_i$'s are D-equivalent each other.
\end{thm}
\noindent
Our Main Theorem means that the case (iii) in Theorem \ref{BMK} really occurs, and 
a similar result to Theorem \ref{Oguiso} is true for elliptic surfaces.
  
In contrast to Main Theorem and Theorem \ref{Oguiso}, 
it is predicted that given a smooth projective variety $X$,
the set of isomorphic classes of FM partners of $X$ is finite.
Actually this is known for the $2$-dimensional case 
(\cite{BM01} and \cite{Ka02}).

\paragraph{Notation and conventions.}
All varieties are defined over $\C$ and ``elliptic surface'' always means ``relatively minimal elliptic surface''
in this note. For a set $I$, we denote by $|I|$ the cardinality of $I$.


\section{The proof of Main Theorem}
We need some standard notation and results before giving the proof.
Let $\pi:S\to C$ be an elliptic surface.
For an object $E$ of $D(S)$, we define the fiber degree of $E$ 
\[d(E)=c_1(E)\cdot f, \]
where $f$ is a general fiber of $\pi$. Let us denote by $\lambda_{S/C}$  
the highest common factor of the fiber degrees of objects of $D(S)$. Equivalently,
$\lambda_{S/ C}$ is the smallest number $d$ such that there is a 
holomorphic $d$-section of $\pi$. 
For integers $a>0$ and $i$ with $i$ coprime to $a\lambda_{S/ C}$, by \cite{Br98} there exists a smooth,
2-dimensional component $J_S (a,i)$ of the moduli space of pure dimension one stable sheaves on $S$,
the general point of which represents a rank $a$, degree $i$ stable vector bundle supported on a smooth 
fiber of $\pi$. 
There is a natural morphism $J_S (a,i)\to C$, taking a point representing a sheaf supported on the
fiber $\pi ^{-1}(x)$ of $S$ to the point $x$. This morphism is a minimal elliptic fibration (\cite{Br98}).
Put $J^i(S):=J_S(1,i)$. 
Obviously, $J^0(S)\cong J(S)$, the Jacobian surface associated to $S$, and $J^1(S)\cong S$. 
 
Fix an elliptic surface with a section $\pi:B\to C$. Let $\eta=\Spec k$ be the generic point of $C$,
where $k=k(C)$ is the function field of $C$, and let $\overline{k}$ be the algebraic closure of $k$.
Put $\overline{\eta}=\Spec \overline{k}$. We define the \emph{Weil--Chatelet group} $WC(B)$
by the Galois cohomology $H^1(G, B_{\eta}(\overline{k}))$. Here $G=\Gal(\overline{k}/k)$ and $B_{\eta}(\overline{k})$ 
is the group of points of the elliptic curve $B_{\eta}$ defined over $\overline{k}$. 
Suppose that we are given a pair $(S,\varphi)$, where $S$ is an elliptic surface $S\to C$ and 
$\varphi$ is an isomorphism $J(S)\to B$ over $C$, fixing their $0$-sections. 
Then we have a morphism 
$$
B_{\eta}\times S_{\eta}\to J(S)_{\eta}\times  S_{\eta}\to S_{\eta}.
$$
Here the first morphism is induced by $\varphi^{-1}\times id_S$ and the second is given by translation.
We obtain a principal homogeneous space $S_{\eta}$ of $B_{\eta}$.
Since this correspondence is invertible and 
the group $H^1(G, B_{\eta}(\overline{k}))$ classifies isomorphic classes 
of principal homogeneous spaces of $B_{\eta}$,
we know that $WC(B)$ consists of all isomorphic classes of pairs $(S,\varphi)$.
Here two pairs $(S,\varphi)$ and $(S',\varphi ')$ are \emph{isomorphic} if there is an isomorphism 
$\alpha :S\to S'$ over $C$, such that $\varphi '\circ \alpha _* =\varphi$,
where $\alpha _*:J(S)\to J(S')$ is the isomorphism induced by $\alpha$ (fixing $0$-sections). 

\noindent
\[ \xymatrix{ J(S) \ar[d]_{\varphi} \ar[r] ^{\alpha _*} & J(S') \ar[d]^{\varphi '} \\
  B \ar@{=}[r] & B}\]

\noindent
There is a short exact sequence (page 185, \cite{Fri98} or page 38, \cite{Fri95}) 
$$
0\to \cy{Sh}(B)\to WC(B) \to \bigoplus_{t \in C} H_{1}(B_{t}, \Q / \Z)\to 0,
$$ 
if $B$ is not the product $C\times E$, where $E$ is an elliptic curve. 
The group $\cy{Sh}(B)$ is called \emph{Tate--Shafarevich group} and it is the subgroup of $WC(B)$ 
which consists of all isomorphic classes of pairs $(S,\varphi)$ such that $S$ does not have multiple fibers.
For a rational surface $B$, it is known that $\cy{Sh}(B)$ is trivial (Example 1.5.12, \cite{FM94}).  

Now we are in position to prove Main Theorem.
\begin{proof}
(i) By the Persson's list \cite{Pe90}, there is a rational elliptic surface $B\to C$ having
a section and three singular fibers of type $III^{*}$, $I_{2}$, $I_{1}$ 
(there are many other choices for $B$). Fix a point $t_{0} \in C$ 
such that $B_{t_{0}}$ is smooth. Take an element $\xi = (\xi_{t})$ of
$WC(B) \cong \bigoplus_{t \in C} H_{1}(B_{t}, \Q / \Z)$
such that $\xi_{t_{0}}$ is of order $p$ and $\xi_{t} = 0$ for other
$t$. Then the surface $\pi:S(p)\to C$ corresponding to $\xi$ is an elliptic surface
with desired singular fibers. We can check that $S(p)$ is rational, for instance,  
by Proposition 1.3.23, \cite{FM94}.

(ii) Put $S=S(p)$. Because every $(-1)$-curve on $S$ is a $p$-section of $\pi$, we know that $\lambda _{S/C}=p$. 
For $i\in \Z$, there is an isomorphism $\varphi_i: J(J^i(S))\to B$ such that 
$(J^i(S),\varphi_i)$ corresponds to $i\xi\in WC(B)$ (\cite{Fri95}, page 38).
By Theorem \ref{BMelliptic}, each $J^i(S)$ is mutually D-equivalent for $1\le i <p$. 
We can also conclude that $J^i(S)$ is rational, since $\kappa (J^i (S))=-\infty$
and the Euler numbers $e(J^i(S))$ and $e(S)$ coincide by Proposition 2.3, \cite{BM01} 
(we can check the rationality also by using Proposition 1.3.23, \cite{FM94}).
Put 
$$
I=\bigl\{1,\ldots,p-1 \bigr\},\qquad
I(a)=\bigl\{i\in I \bigm| J^i(S)\cong J^a(S)\bigr\}
$$ 
for $a\in I$. Then there are $i_1,\ldots,i_M\in I$ such that $I=\coprod_{k=1}^{M}I(i_k)$ (disjoint union).
\begin{cla}\label{claim}
For all $a\in I$, $|I(a)|\le 6$.
\end{cla}
\noindent
If Claim \ref{claim} is true, we have $6M\ge |I|=p-1$.
By the assumption $p>6(N-1)+1$, we have $M\ge N$, which completes the proof of Main Theorem. 

Let us start the proof of Claim \ref{claim}.
\paragraph{Step 1.}
For each $i\in I(a)$, we fix an isomorphism $\alpha_i:J^a(S)\to J^i(S)$. 
Because the rational surface $J^a(S)$ has a unique elliptic fibration,
there exists $\delta \in\Aut C$ such that the following diagram is commutative.

\noindent
\[ \xymatrix{ J^a(S) \ar[d] \ar[r]^{\alpha_i}& J^i (S) \ar[d] \\
  C \ar[r] ^{\delta} & C}\]

\noindent
This makes the following diagram commutative.

\noindent
\[ \xymatrix{J(J^a(S))\  \ar[d] \ar[r]^{\varphi _a ^{-1}\circ\varphi _i\circ\alpha _{i*}} & \ J(J^a (S)) \ar[d] \\
  C \ar[r] ^{\delta} & C}\]

\noindent
By our assumption, $J(J^a(S))$ has at least three singular fibers of different Kodaira's types.
Hence $\delta$ must be the identity on $C\cong\PP ^1$ and 
then we can say that every $\alpha_i$ is an isomorphism over $C$.    

\paragraph{Step 2.} 
By Step 1, we know that $\varphi_i\circ\alpha _{i *}\circ\varphi^{-1}_a$ is an automorphism of $B$ over $C$,
fixing the $0$-section.
Put $\gamma_i=\varphi_i\circ\alpha _{i *}\circ\varphi^{-1}_a$.

\noindent
\[ \xymatrix{ J(J^a(S)) \ar[d]_{\varphi _a} \ar[r] ^{{\alpha _i}_*} & J(J^i (S)) \ar[d]^{\varphi _i} \\
  B \ar[r] ^{\gamma _i} & B}\]

\noindent
Suppose $\gamma_i=\gamma_j$ for $i,j\in I(a)$,
then by the isomorphism $\alpha_j\circ\alpha_i^{-1}$, 
we see that $(J^i(S),\varphi_i)$ is isomorphic to $(J^j(S),\varphi_j)$ and hence   
$i\xi=j\xi$ in $WC(B)$. 
Because the order of $\xi$ is $p$, we obtain $i=j$. 
Since the order of the group of
automorphism of $B$ over $C$ fixing the 0-section is at most 6, 
we get $|I(a)|\le 6$. This finishes the proof.
\end{proof}

\begin{thm}[Proposition 4.4, \cite{BM01}]\label{BMelliptic}
Let $\pi :S\to C$ be an elliptic surface and $T$ a smooth projective variety.
Assume that the Kodaira dimension $\kappa (S)$ is non-zero.
Then the following are equivalent.
\renewcommand{\labelenumi}{(\roman{enumi})}
\begin{enumerate}
\item 
$T$ is an FM partner of $S$. 
\item 
$T$ is isomorphic to $J^b(S)$ for some integer $b$ with $(b,\lambda _{S/C})=1$. 
\end{enumerate}
\end{thm}

 I would like to express my hearty thanks to Professors Akira Ishii, Noboru Nakayama and Takeshi Kajiwara
for useful conversation. 
I am very grateful to the referee for his/her invaluable suggestions.

\indent

\textsc{Department of Mathematics, Kyoto University, Kyoto, 606-8502, Japan} \\
\texttt{hokuto@kurims.kyoto-u.ac.jp}


\begin{thebibliography}{Fri95}

\bibitem{Br98}
T. Bridgeland, Fourier--Mukai transforms for elliptic surfaces, J. Reine Angew. Math. 498 (1998), 115--133. 
\bibitem{BM01}
T. Bridgeland, A. Maciocia, Complex surfaces with equivalent derived categories, Math. Z. 236 (2001), 677--697. 
\bibitem{Fri95}
R. Friedman, Vector bundles and $SO(3)$-invariants for elliptic surfaces, J. Amer. Math. Soc. 8, (1995), 29-139.
\bibitem{Fri98}
R. Friedman, Algebraic surfaces and holomorphic vector bundles. Universitext. Springer-Verlag, New York, 1998.
\bibitem{FM94}
R. Friedman, J. Morgan, Smooth four-manifolds and complex surfaces. Ergebnisse der Mathematik und ihrer Grenzgebiete (3),
27. Springer-Verlag, Berlin, 1994. x+520 pp
\bibitem{HLOY}
S. Hosono, B. H. Lian, K. Oguiso, S.-T. Yau. 
Kummer structures on a K3 surface - An old question of T. Shioda, math.AG/0202082
\bibitem{Ka02}
Y. Kawamata, D-equivalence and K-equivalence, J. Differential Geom. 61 (2002), 147-171.

\bibitem{Oguiso02}
K. Oguiso, K3 surfaces via almost-primes, Math. Res. Lett. 9 (2002), 47--63.
\bibitem{Pe90}
U. Persson, Configurations of Kodaira fibers on rational elliptic surfaces. Math. Z. 205 (1990), 1--47.
\end{thebibliography}
\end{document}